\documentclass[reqno, 12pt]{amsart}
\usepackage{amssymb,amsmath,amsfonts,amsthm,comment,mathrsfs,times,graphicx}
\usepackage[bookmarksnumbered, colorlinks, plainpages]{hyperref}
\usepackage{color}
\usepackage[english]{babel}
\usepackage[all,cmtip]{xy}
\usepackage{lmodern}
\usepackage{enumitem}
\usepackage{geometry}
\newcounter{dummy} \numberwithin{dummy}{section}
\newtheorem{theo}[dummy]{Theorem }
 \newtheorem{coro}[dummy]{Corollary}
 
 \newtheorem{pro}[dummy]{Proposition}
\newtheorem{deft}[dummy]{Definition}

\newtheorem{rem}[dummy]{Remark}

\newtheorem{conj}{Conjecture}
\title[ Fibonacci sequences and real quadratic $p$-rational fields]{Fibonacci sequences and real quadratic $p$-rational fields}
\author[ Zakariae Bouazzaoui ]{ Zakariae Bouazzaoui$^{(1)}$}
\address{$^{(1)}$ Université Moulay Ismaïl,
 Département de mathématiques,
Faculté des sciences de Meknès, B.P. 11201 Zitoune, Meknès, Maroc.}
\email{\textcolor[rgb]{0.00,0.00,1.00}{z.bouazzaoui@edu.umi.ac.ma}}

 
\begin{document}
\maketitle
\renewcommand{\abstractname}{Abstract}
\begin{abstract}
We study the $p$-rationality of real quadratic fields in terms of
generalized Fibonacci numbers and their periods modulo positive
integers.
\end{abstract}

\section{Introduction}

Let $K$ be a number field and $p$ an odd prime number. The field $K$
is said to be $p$-rational if the Galois group of the maximal
pro-$p$-extension of $K$ which is unramified outside $p$ is a free
pro-$p$-group of rank $r_2+1$, where $r_2$ is the number of pairs of
complex embedding of $K$. The notion of $p$-rational number fields
has been introduced by Movahhedi and Nguyen Quang Do
\cite{Movahhedi-Nguyen}, \cite{Movahhedi}, \cite{Movahhedi90}, and
is used for the construction of non-abelian extensions satisfying
Leopoldt's conjecture. Recently, R. Greenberg used complex abelian
$p$-rational number fields for the construction of $p$-adic Galois
representations with open images. In these paper we study the
$p$-rationality of real quadratic number fields. In fact, we give a
generalization of a result of Greenberg \cite[Corollary
4.1.5]{Greenberg} which relates the $p$-rationality of the field
$\mathbf{Q}(\sqrt{5})$ to properties of the classical Fibonacci
numbers. More precisely, let $d>0$ be a fundamental discriminant.
Denote by $\varepsilon_d$ and $h_d$ the fundamental unit and the
class number of the field $\mathbf{Q}(\sqrt{d})$ and let $N(.)$ be
the absolute norm. We associate to the field $\mathbf{Q}(\sqrt{d})$
a Fibonacci sequence
$F^{(\varepsilon_d+\overline{\varepsilon}_d,N(\varepsilon_d))}=(F_n)_{n\geq0}$
defined by $F_0=0$, $F_1=1$ and the recursion formula
$$F_{n+2}=(\varepsilon_d+\overline{\varepsilon}_d)F_{n+1}-N(\varepsilon_d)F_{n},\; for\;n\geq0.$$

The main result of this paper is the following theorem, which
describes the $p$-rationality in terms of Fibonacci-Wieferich prime
(see Definition \ref{fibonacci-wieferich} for Fibonacci-Wieferich
primes).

\begin{theo}\label{Main theorem}
Let $p\geq3$ be an odd prime number such that
$p\nmid(\varepsilon_d-\overline{\varepsilon}_d)^{2}h_d$. The
following assertions are equivalent:
\begin{enumerate}
    \item  the field $\mathbf{Q}(\sqrt{d})$ is $p$-rational,
    \item $p$ is not a Fibonacci-Wieferich prime for
    $\mathbf{Q}(\sqrt{d})$.
\end{enumerate}
\end{theo}

It is known that for every positive integer $m$, the reduction
modulo $m$ of the sequence $(F_{n})_n$ is periodic of period a
positive integer $k(m)$ \cite[Theorem 1.]{Wall}, \cite{Dan-Robert}.
Using this fact and properties of these periods, we give another
characterization of the $p$-rationality of $\mathbf{Q}(\sqrt{d})$ in
terms of the periods of the associated Fibonacci numbers.

\begin{pro}\label{periods and p-rationality}
Let $p\geq3$ be a prime number such that $p\nmid
(\varepsilon_{d}-\overline{\varepsilon_{d}})^2 h_d$. Then the
following assertions are equivalent:
\begin{enumerate}
    \item  the field $\mathbf{Q}(\sqrt{d})$ is $p$-rational,
    \item $k(p)\neq k(p^2)$.
\end{enumerate}
\end{pro}

For the classical Fibonacci sequence, namely $a=b=1$, D.D. Wall is
the first to study these periods in \cite{Wall}, where he proved
many properties of these integers. One problem encountered by Wall
in his paper is the study of the hypothesis $k(p)\neq k(p^2)$. He
asked whether the equality $k(p)=k(p^2)$ is possible. This question
is still open with strong numerical evidence
\cite{Elsenhans-Jahnel}.
By Proposition \ref{periods and p-rationality}, it is equivalent to
whether the number field $\mathbf{Q}(\sqrt{5})$ is not $p$-rational
for some prime number $p$. It is generalized to Fibonacci sequences
$F^{(a,b)}$ where for some sequences we have an affirmative answer,
for example the Fibonacci sequence $F^{(2,-1)}$ gives that
$k(13)=k(13^2)$ and $k(31)=k(31^2)$, which means that the field
$\mathbf{Q}(\sqrt{2})$ is not $p$-rational for $p=13,31$. Under the
light of the above characterization of the $p$-rationality, the
conjecture of G. Gras \cite[Conjecture 7.9]{Gras} on the
$p$-rationality of real quadratic fields, it is suggested that for
almost all primes $p$ we have $k(p)\neq k(p^2)$.\\

\section{$p$-rational fields}

In this section we give a characterization of the $p$-rationality of
real quadratic fields in terms of values of the associated
$L$-functions at odd negative integers. In fact, the $p$-rationality
of totally real abelian number fields $K$ is intimately related to
special values of the associated zeta functions $\zeta_K$. The
relation is as follows. For any finite set $\Sigma$ of primes of
$K$, we denote by $G_{\Sigma}(K)$ the Galois group of the maximal
pro-$p$-extension of $K$ which is unramified outside $\Sigma$. Let
$S$ be the finite set of primes $S_{p}\cup S_{\infty}$, where
$S_{\infty}$ is the set of infinite primes of $K$ and $S_{p}$ is the
primes above $p$ in $K$. It is known that the group $G_{S_p}(K)$ is
a free pro-$p$-group on $r_2+1$ generators if and only if the second
Galois cohomology group $H^2(G_{S_p}(K),\mathbf{Z}/p\mathbf{Z})$
vanishes. This vanishing is related to special values of the zeta
function $\zeta_K$ via the conjecture of Lichtenbaum. More
precisely, let $\mathcal{G}_S$ be the Galois group of the maximal
extension of $K$ which is unramified outside $S$. The main
conjecture of Iwasawa theory (now a theorem of Wiles \cite{Wiles90})
relates the order of the group $H^2(\mathcal{G}_S,\mathbf{Z}_p(i))$,
for even integers $i$, to the $p$-adic valuation of $\zeta_K(1-i)$
by the $p$-adic equivalence:
\begin{equation}\label{Main conjecture}
w_{i}(K)\zeta_K(1-i)\sim_p |H^2(\mathcal{G}_S,\mathbf{Z}_{p}(i))|,
\end{equation}
where for any integer $i$, $w_i(F)$ is the order of the group
$H^0(G_F,\mathbf{Q}_p/\mathbf{Z}_p(i))$, and $\sim_p$ means having
the same $p$-adic valuation, see e.g \cite{Kolster}. Moreover, the
group $H^2(\mathcal{G}_S,\mathbf{Z}_{p}(i))$ vanishes if and only if
$H^2(\mathcal{G}_S,\mathbf{Z}/p\mathbf{Z}(i))$ vanishes. Let $\mu_p$
be the group of $p$-th unity. The periodicity of the groups
$H^2(\mathcal{G}_S,\mathbf{Z}/p\mathbf{Z}(i))$ modulo
$\delta=[K(\mu_p):K]$ gives that
$$H^2(\mathcal{G}_S,\mathbf{Z}/p\mathbf{Z}(i))\cong
H^2(\mathcal{G}_S,\mathbf{Z}/p\mathbf{Z}(i+j\delta)),$$ for any
integer $j$. In addition, since $p$ is odd, the vanishing of the
group $H^2(\mathcal{G}_S,\mathbf{Z}/p\mathbf{Z}(i))$ is equivalent
to the vanishing of the group
$H^2(G_{S_p}(K),\mathbf{Z}/p\mathbf{Z}(i))$. Number fields such that
$H^2(G_{S_p}(K),\mathbf{Z}/p\mathbf{Z}(i))=0$ are called
$(p,i)$-regular \cite{Assim}. In particular, the field $K$ is
$p$-rational if and only if $w_{p-1}(K)\zeta_K(2-p)\sim_{p} 1$. This
leads to the following characterization of the $p$-rationality of
totaly real number fields.

\begin{pro}\label{p-rational and L-fun}
Let $p$ be an odd prime number which is unramified in an abelian
totally real number field $K$. Then we have the equivalence
\begin{equation}
K\;\hbox{is p-rational}\;\Leftrightarrow\; L(2-p,\chi)\;\hbox{is a
p-adic unit},
\end{equation}
where $\chi$ is ranging over the set of irreducible characters of
$\mathrm{Gal}(K/\mathbf{Q})$.
\end{pro}

\noindent{\textbf{Proof}.} First, the zeta function $\zeta_{K}$
decomposes in the following way:
$$\zeta_{K}(2-p)=\zeta_{\mathbf{Q}}(2-p)\times\prod_{\chi\neq1}L(2-p,\chi).$$
Second, it is known that $\zeta_{\mathbf{Q}}(2-p)$ is of $p$-adic
valuation $-1$ and that $w_{p-1}(K)$ has $p$-adic valuation $1$,
giving that $w_{p-1}(K)\zeta_{\mathbf{Q}}(2-p)\sim_p 1$. Then from
(\ref{Main conjecture}) we obtain the formula
$$\prod_{\chi\neq1}L(2-p,\chi)\sim_p |H^2(\mathcal{G}_S,\mathbf{Z}_{p}(p-1))|.$$
Since, for every character $\chi\neq1$, the value $L(2-p,\chi)$ is a
$p$-integers \cite[Corollary 5.13]{Washington}, we have
$H^2(\mathcal{G}_S,\mathbf{Z}_p(p-1))=0$ if and only if for every
$\chi\neq1$, $L(2-p,\chi)$ is a $p$-adic unit. Furthermore, the
vanishing of the group $H^2(\mathcal{G}_S,\mathbf{Z}_{p}(p-1))$ is
equivalent to the vanishing of the group
$H^2(\mathcal{G}_S,\mathbf{Z}/p\mathbf{Z}(p-1))$, which turns out to
be equivalent to the vanishing of
$H^2(G_{S_p}(K),\mathbf{Z}/p\mathbf{Z})$ (by the above mentioned
periodicity statement). This last vanishing occurs exactly when the
field $K$ is $p$-rational.\hfill $\blacksquare$\vskip 6pt

In the particular case of a real quadratic field
$K=\mathbf{Q}(\sqrt{d})$, we have the decomposition
$$\zeta_{K}(2-p)=\zeta_{\mathbf{Q}}(2-p)L(2-p,(\frac{d}{.})),$$
where $(\frac{d}{.})$ is the quadratic character associated to the
field $K=\mathbf{Q}(\sqrt{d})$.

\begin{coro}\label{L-function}
For every odd prime number $p\nmid d$, we have the equivalence
\begin{equation}
\mathbf{Q}(\sqrt{d})\;\hbox{is
p-rational}\;\Leftrightarrow\;L(2-p,(\frac{d}{.}))\not\equiv0\pmod{p}.
\end{equation}
\end{coro}\hfill $\blacksquare$\vskip 6pt

\begin{rem}\label{remark 1}
The properties of special values of $p$-adic $L$-functions tells us
that the $p$-rationality is related to the class number and the
$p$-adic regulator. More precisely, let $K$ be a totally real number
field of degree $g$. Under the Leopoldt conjecture, class field
theory gives that
$G_{S_p}(K)^{ab}\cong\mathbf{Z}_{p}^{r_2+1}\times\mathcal{T}_K$,
where $\mathcal{T}_K$ is the $\mathbf{Z}_p$-torsion of
$G_{S_p}(K)^{ab}$. Then the field $K$ is $p$-rational precisely when
$\mathcal{T}_K=0$ \cite[Théor\`{e}me et Definition
1.2]{Movahhedi-Nguyen}. Moreover, the order of $\mathcal{T}_K$
satisfies
\begin{equation}
|\mathcal{T}_K|\sim_p
w(K(\mu_p))\prod_{v|p}(1-N(v)^{-1}).\frac{R_p(K).h_K}{\sqrt{|d_k|}},
\end{equation}
(\cite[app]{Coates}), where $h_K$ is the class number, $R_p(K)$ is
the $p$-adic regulator, $N(v)$ is the absolute norm of $v$,
$w(K(\mu_p))=|\mu(K(\mu_p))|$ the number of roots of unity of
$K(\mu_p)$ and $d_K$ is the discriminant of the number field $K$.
Hence for every odd prime number $p$ such that $(p,d_{K}h_K)=1$, the
field $K$ fails to be $p$-rational if and only if $v_p(R_p(K))>g-1$.
\end{rem}

Under the light of Remark \ref{remark 1}, for a real quadratic field
$\mathbf{Q}(\sqrt{d})$ we have the equivalence
\begin{equation}\label{regulator and p-rationality}
\mathbf{Q}(\sqrt{d})\;\hbox{is
p-rational}\;\Leftrightarrow\;R_p(\mathbf{Q}(\sqrt{d}))\not\equiv0\pmod{p^{2}}.
\end{equation}
Recall that $R_p(\mathbf{Q}(\sqrt{d}))=\log_{p}{(\varepsilon_d)}$,
where $\varepsilon_d$ is a fundamental unit of $K$ and $\log_p$ is
the $p$-adic logarithm.

\section{Fibonacci number}

The classical Fibonacci sequence is an interesting linear recurrence
sequence, in part because of its applications in several areas of
sciences. Here we consider a class of linear recurrence sequences
which arise from real quadratic fields and that we use for the study
of the $p$-rationality of these fields. As mentioned in the
introduction, Greenberg \cite[Corollary 4.1.5.]{Greenberg} used
classical Fibonacci numbers to give a characterization for the
$p$-rationality of the field $\mathbf{Q}(\sqrt{5})$. In this paper
we give a generalization of this result to any real quadratic field.
The Fibonacci numbers associated to real quadratic fields are given
as follows. Let $d>0$ be a fundamental discriminant and let $h_{d}$,
$\varepsilon_d$ be respectively the class number and the fundamental
unit of the field $\mathbf{Q}(\sqrt{d})$ with ring of integers
$\mathcal{O}_d$. We denote by $\overline{\varepsilon}_d$ the
conjugate of $\varepsilon_d$ and $N(.)$ the absolute norm. Define
the sequence
$F^{(\varepsilon_d+\overline{\varepsilon}_d,N(\varepsilon_d))}=(F_n)_n$
such that $F_0=0$, $F_1=1$ and
$$F_{n+2}=(\varepsilon_d+\overline{\varepsilon}_d)F_{n+1}-N(\varepsilon_d)F_n.$$
The Binet formula \cite[page 173]{Dan-Robert} gives that
$$F_n=\frac{\varepsilon_d^{n}-\overline{\varepsilon}_d^{n}}{\varepsilon_d-\overline{\varepsilon}_d},\quad\quad \forall n\geq0.$$
We establish a relation between Fibonacci numbers and the $p$-adic
regulator which allows us to prove the main result.

\begin{deft}
Let $a$ be a non trivial element of the ring of integers of the
field $\mathbf{Q}(\sqrt{d})$ such that $(a,p)=1$. Then the prime $p$
is said to be Wieferich of basis $a$ if the following congruence
holds:
$$a^{p^{r}-1}-1\equiv0\pmod{p^{2}},$$ where $r$ is the residue degree of $p$
in the quadratic field $\mathbf{Q}(\sqrt{d})$. Otherwise, the prime
number $p$ is said to be non-Wieferich of basis $a$.
\end{deft}

We have the following equality
$$\log_p{((\varepsilon_{d}^{p^{r}-1}-1)+1)}=(\varepsilon_{d}^{p^{r}-1}-1)-\frac{1}{2}(\varepsilon_{d}^{p^{r}-1}-1)^2+...$$
where $\log_p$ is the $p$-adic logarithm and as before $r$ is the
residue degree of $p$ in the quadratic field $\mathbf{Q}(\sqrt{d})$.
Since $R_p=\log_p{(\varepsilon_{d})}$ and the group
$(\mathcal{O}_d/p\mathcal{O}_d)^{\times}$ is cyclic of order
$p^r-1$, where $\mathcal{O}_d$ is the ring of integers of
$\mathbf{Q}(\sqrt{d})$, we obtain the equivalences
\begin{equation}\label{Wieferich-regulator}
  \begin{array}{ccc}
    \varepsilon_d^{p^{r}-1}-1\not\equiv0 \pmod{p^{2}} & \Leftrightarrow & R_p\equiv p\pmod{p^{2}}, \\
     & \Leftrightarrow & R_p\not\equiv0\pmod{p^{2}}. \\
  \end{array}
\end{equation}

Then combining this last equivalence with the equivalence
(\ref{regulator and p-rationality}) we obtain

\begin{pro}
Let $p$ be an odd prime number such that $p\nmid dh_d$. Then the
field $\mathbf{Q}(\sqrt{d})$ is $p$-rational if and only if $p$ is a
non-Wieferich prime of basis $\varepsilon_d$.
\end{pro}\hfill $\blacksquare$\vskip 6pt

Very little is known about these primes and it is conjectured that
the set of Wieferich primes is of density zero \cite{Silverman}. In
the following we are interested with the set of Fibonacci-Wieferich
primes defined as follows.

\begin{deft}\label{fibonacci-wieferich}
A prime number $p$ is said to be a Fibonacci-Wieferich prime for the
field $\mathbf{Q}(\sqrt{d})$ if
$$F_{p-(\frac{d}{p})}\equiv0\pmod{p^2},$$
where $(\frac{d}{.})$ is the Legendre symbol associated to the
quadratic field $\mathbf{Q}(\sqrt{d})$.
\end{deft}

We give the main result of this section which describe the
$p$-rationality in terms of Fibonacci-Wieferich primes.

\begin{theo}\label{fibonacci}
Let $p\geq5$ be a prime number such that $p\nmid
(\varepsilon_{d}-\overline{\varepsilon_{d}})^2 h_d$. Then the
following assertions are equivalent:
\begin{enumerate}
    \item  the field $\mathbf{Q}(\sqrt{d})$ is $p$-rational,
    \item $p$ is not a Fibonacci-Wieferich prime for
    $\mathbf{Q}(\sqrt{d})$.
\end{enumerate}
\end{theo}

\noindent{\textbf{Proof.}} Using the equivalence
(\ref{Wieferich-regulator}), it suffices to prove that:
\begin{equation}
\varepsilon_d^{p^{r}-1}-1\not\equiv0\pmod{p^{2}}\;\Leftrightarrow\;F_{p-(\frac{d}{p})}\not\equiv0\pmod{p^{2}}.
\end{equation}
Let $Q_p(\varepsilon_d)$ be the residue class
$$\frac{\varepsilon_d^{p^{r}-1}-1}{p}\pmod{p}.$$
A prime number $p$ satisfying
$Q_p(\varepsilon_d)\not\equiv0\pmod{p}$ is non-Wieferich of basis
$\varepsilon_d$.\\
First suppose that $(\frac{d}{p})=1$. Then $r=1$ and
\begin{equation}
Q_p(\varepsilon_d)\equiv \frac{\varepsilon_d^{p-1}-1}{p}\pmod{p}.
\end{equation}
The Binet formula gives that
$$(\varepsilon_d-\overline{\varepsilon}_d)F_{p-1}=\varepsilon_d^{p-1}-\overline{\varepsilon}_d^{p-1}=
\varepsilon_d^{1-p}(\varepsilon_d^{(p-1)}-1)(\varepsilon_d^{(p-1)}+1).$$
Since $\varepsilon_d$ is a unit and $p\nmid
(\varepsilon_d-\overline{\varepsilon}_d)$, we have
$\varepsilon_d^{1-p}(\varepsilon_d^{(p-1)}+1)(\varepsilon_d-\overline{\varepsilon}_d)\not\equiv0\pmod{p}$.
Hence we obtain the equivalence
\begin{equation}
Q_p(\varepsilon_d)\not\equiv0\pmod{p}\Leftrightarrow
F_{p-1}\not\equiv0\pmod{p^2}.
\end{equation}
Second, suppose that the prime number $p$ is inert in the field
$\mathbf{Q}(\sqrt{d})$. Then we have
\begin{equation}
Q_p(\varepsilon_d)\equiv \frac{\varepsilon_d^{p^2-1}-1}{p}\pmod{p}.
\end{equation}
The Galois group of $\mathbf{Q}(\sqrt{d})/\mathbf{Q}$ is generated
by an element $\sigma$ of order two such that
$\sigma(\varepsilon_d)=\overline{\varepsilon}_d$. Since the group
$(\mathcal{O}_d/p\mathcal{O}_d)^{\times}$ is cyclic of order
$p^2-1$, we have $\varepsilon_d^{p+1}\equiv x\pmod{p}$ for some
$x\in\mathbf{Z}$. Hence $\overline{\varepsilon}_d^{p+1}\equiv
x\pmod{p}$ and $F_{p+1}\equiv0\pmod{p}$. Note that since
$$F_{p+1}=(\varepsilon_d-\overline{\varepsilon}_d)^{-1}\overline{\varepsilon}_d^{p+1}(\varepsilon_d^{2(p+1)}-1),$$
we have $$\varepsilon_d^{2(p+1)}-1\equiv0\pmod{p}.$$ Moreover,
$$Q_p(\varepsilon_d)=\frac{1}{p}(\varepsilon_d^{p^{2}-1}-1)=\frac{1}{p}((\varepsilon_d^{2(p-1)})^{\frac{p-3}{2}}-1)=
\frac{1}{p}(\varepsilon_d^{2(p+1)}-1)(\varepsilon_d^{2(p+1)\frac{p-1}{2}}+...+1).$$
Since $x^{2}\equiv1\pmod{p}$, we then obtain the congruence
$$Q_p(\varepsilon_d)\equiv\frac{1}{p}\frac{p-1}{2}(\varepsilon_d^{2(p+1)}-1)\pmod{p}.$$
Hence we have the equivalence
\begin{equation}
Q_p(\varepsilon_d)\not\equiv0\pmod{p}\Leftrightarrow
F_{p+1}\not\equiv0\pmod{p^2}.
\end{equation}
Then in all cases we obtain that the field $\mathbf{Q}(\sqrt{d})$ is
$p$-rational precisely when $p$ is not a Fibonacci-Wieferich prime.
\hfill $\blacksquare$\vskip 6pt

Using this characterization of the $p$-rationality on pariGP, we
obtain some numerical evidence for the primes $p$ for which a given
real quadratic number field is not $p$-rational.
$$
\begin{tabular}{|c|c|}
  \hline
  Discriminant & Primes$<10^9$\\
  \hline
5&  \\
\hline
8& 13, 31, 1546463 \\
\hline
12& 103 \\
\hline
13& 241 \\
\hline
17&  \\
\hline
21& 46179311 \\
\hline
24& 7, 523 \\
\hline
28&  \\
\hline
29& 3, 11 \\
\hline
33& 29, 37, 6713797 \\
\hline
37& 7, 89, 257, 631 \\
\hline
40& 191, 643, 134339, 25233137 \\
\hline
41& 29, 53, 7211 \\
\hline
44&  \\
\hline
53& 5 \\
\hline
56& 6707879, 93140353 \\
\hline
57& 59, 28927, 1726079, 7480159 \\
\hline
60& 181, 1039, 2917, 2401457 \\
\hline
61&  \\
\hline
65& 1327, 8831, 569831 \\
\hline
69& 5, 17, 52469057 \\
\hline
73& 5, 7, 41, 3947, 6079 \\
\hline
76& 79, 1271731, 13599893, 31352389\\
\hline
77& 3, 418270987 \\
\hline
85& 3, 204520559 \\
\hline
88& 73, 409, 43, 28477 \\
\hline
89& 5, 7, 13, 59 \\
\hline
92& 7, 733 \\
\hline
93& 13 \\
\hline
97& 17, 3331\\
\hline
\end{tabular}
$$\\

With the help of these results and further computations, we could
construct examples of multi-quadratic $p$-rational fields. The first
example is the field
$K_1=\mathbf{Q}(\sqrt{2},\sqrt{3},\sqrt{5},\sqrt{7},\sqrt{11},\sqrt{-1})$,
which is $p$-rational for all primes $100<p<1000$ except for
$p=103,173,181,191,199,227,251,269,\\409,523,571,577,643,859$.
Another example is the field
$K_2=\mathbf{Q}(\sqrt{13},\sqrt{17},\sqrt{19},\sqrt{23},\sqrt{29},\sqrt{-1})$.
The field $K_2$ is $p$-rational for all primes $100<p<1000$ except
for
$151,197,227,241,307,337,401,457,\\487,593,643,709,719,733,809,839$.
Hence for every prime $100<p<1000$ such that $p\neq227,643$, there
exist a $p$-rational field of degree $2^{t}$ for any $1\leq
t\leq6$.\\
The above examples are weak numerical evidence to a conjecture
proposed by Greenberg:

\begin{conj}(\cite[Conjecture
4.2.1.]{Greenberg}) For any odd prime $p$ and for any $t\geq1$,
there exists a $p$-rational field $K$ such that
$\mathrm{Gal}(K/\mathbf{Q})\cong(\mathbf{Z}/2\mathbf{Z})^{t}$.
\end{conj}

As an important consequence of this conjecture, Greenberg proved the
following proposition.

\begin{pro}\cite[Proposition 6.2.2]{Greenberg}
Suppose that $K$ is a complex $p$-rational number field and that
$\mathrm{Gal}(K/\mathbf{Q})$ is isomorphic to
$(\mathbf{Z}/2\mathbf{Z})^{t}$, where $t\geq4$. Let $n$ be an
integer such that $4\leq n\leq 2^{t-1}-3$. Then there exists a
continuous homomorphism $$\rho: G_{\mathbf{Q}}\rightarrow
GL_n(\mathbf{Z}_p),$$ with open image.
\end{pro}

Based on the above computations and Proposition $3.5$, we have the
following corollary.

\begin{coro}
For any integer $4\leq n\leq 2^{5}-3$ and any prime $100<p<1000$
such that $p\neq227,643$, there exists a $p$-adic Galois
representation
$$\rho: G_{\mathbf{Q}}\rightarrow GL_n(\mathbf{Z}_p),$$ with open
image.
\end{coro}

Another characterization of the $p$-rationality is given in terms of
periods of Fibonacci sequences modulo $p$ and $p^{2}$. Let
$F^{(a,b)}$ be a Fibonacci sequence and $m$ a positive integer such
that $(b,m)=1$. As mentioned above the sequence $F^{(a,b)}\pmod{m}$
is periodic of period $k(m)$. Wall studied these periods for
classical Fibonacci sequence and general results are obtained in
\cite[page 374-376]{Renault}. We describe the $p$-rationality of
real quadratic fields in terms of periods of Fibonacci sequence
associated to these fields.

\begin{theo}\label{periods}\cite[Proposition 3.2.4]{Elsenhans-Jahnel}
\begin{equation}
\hbox{The equality}\; k(p)=k(p^2)\;\hbox{holds}\; \hbox{if and only
if}\; F_{p-(\frac{d}{p})}\equiv0\pmod{p^2}.
\end{equation}
\end{theo}

Proposition \ref{periods and p-rationality} follows from Theorem
\ref{fibonacci} and Theorem \ref{periods}. For the classical
Fibonacci numbers $F_n$, the field $\mathbf{Q}(\sqrt{5})$ is
$p$-rational precisely when $p$ is not a Fibonacci-Wieferich prime
\cite[Corollary 4.1.5]{Greenberg}. It is known that up to
$6.7\times10^{15}$ there is no Fibonacci-Wieferich primes
\cite{François-Klyve}. Greenberg pointed out in \cite{Greenberg}
that such primes are quite rare, they have trivial density if we
assume G. Gras Conjecture, which asserts that a number field is
$p$-rational for almost all primes. Theorem \ref{periods}, gives
that the field $\mathbf{Q}(\sqrt{d})$ is $p$-rational if and only if
$k(p)\neq k(p^2)$. According to the table above there is fundamental
discriminants $d$ such that there exist primes $p$ for which
$k(p)=k(p^2)$. As an example we mentioned the case of
$\mathbf{Q}(\sqrt{2})$ where $k(13)=k(13^2)$. Note that up to
$10^{9}$, for some discriminants we still have no primes satisfying
the equality of Wall such as $17,28,44,61$.

\section{Williams Congruence}

Let $d$ be a positive fundamental discriminant and $p$ be an odd
prime number such that $p\nmid d$. We are interested with the
numbers $F_{p-(\frac{d}{p})}$. In the classical case, namely the
field $\mathbf{Q}(\sqrt{5})$, we have explicit formula for the
quotient $F_{p-(\frac{5}{p})}/p$ \cite[Theorem 4.1]{Williams}. For
the general case we have a result due to H.C. Williams in
\cite{Williams} which describes these quotients for any real
quadratic field. The results obtained in the above section, combined
with the formula proved by Williams gives another characterization
of the $p$-rationality of real quadratic fields. For an integer $n$,
let $\{n\}$ be the least non-negative residue of $n$ modulo $d$. The
integer $p'$ represents the inverse of $p$ modulo $d$ and
$(\frac{d}{.})$ is the Legendre symbol. Consider the following sum
of characters:
\begin{equation*}
\beta_p(i)=\sum_{j=1}^{\{p'i\}-1}(\frac{d}{j}).
\end{equation*}
Then the result of Williams is as follows:

\begin{theo}\label{Williams}\cite{Williams} Let $p$ be an odd prime number such
that $p\nmid d$. Then
\begin{equation}\label{Williams congruence}
h_dF_{p-(\frac{d}{p})}/p\equiv-2(\frac{d}{p})N^{\frac{(\frac{d}{p})-1}{2}}
\sum_{i=1}^{\frac{p-1}{2}}\beta_p(i)\frac{1}{i}\pmod{p},
\end{equation}
where $h_d$ is the class number of the field $\mathbf{Q}(\sqrt{d})$,
and $\frac{1}{i}$ is the inverse of $i$ modulo $p$.
\end{theo}

An interesting problem of combinatorics and additive number theory
is the study of sums of reciprocals in finite fields. Here we are
concerned with the linear combinations
$$\sum_{i=1}^{d}\beta_p(i)\alpha_p(i)\pmod{p},$$
where $$\alpha_p(i)=\sum_{\begin{array}{c}
                      1\leq k\leq \frac{p-1}{2} \\
                      k\equiv i\pmod{d}
                    \end{array}
}\frac{1}{k}.$$ We have the following description of the
$p$-rationality of the field $\mathbf{Q}(\sqrt{d})$.

\begin{theo}\label{sum of reciprocals}
If $p$ does not divides
$(\varepsilon_d-\overline{\varepsilon}_d)^2$, then
\begin{equation}
\mathbf{Q}(\sqrt{d})\;\hbox{is
p-rational}\;\Leftrightarrow\;\sum_{i=1}^{d}\beta_p(i)\alpha_p(i)\not\equiv0\pmod{p}.
\end{equation}
\end{theo}

\noindent{\textbf{Proof}.} It is known that
$F_{p-(\frac{d}{p})}\equiv0\pmod{p}$ \cite[page 431 formula
(1.2)]{Williams}. Then by Theorem \ref{Main theorem}, the field
$\mathbf{Q}(\sqrt{d})$ is $p$-rational if and only if
$h_{d}F_{p-(\frac{d}{p})}/p\not\equiv0\pmod{p}$. Using Theorem
\ref{Williams}, this occurs precisely when
$$-2(\frac{d}{p})N^{\frac{(\frac{d}{p})-1}{2}}
\sum_{i=1}^{\frac{p-1}{2}}\beta_p(i)\frac{1}{i}\not\equiv0\pmod{p}.$$
Since $p$ is an odd prime number and the term
$2(\frac{d}{p})N^{\frac{(\frac{d}{p})-1}{2}}$ equals $1$ or $2$, the
field $\mathbf{Q}(\sqrt{d})$ is $p$-rational if and only if
\begin{equation}\label{Main equivalence}
\sum_{i=1}^{\frac{p-1}{2}}\beta_p(i)\frac{1}{i}\not\equiv0\pmod{p}.
\end{equation}
Recall that for any integer $i$, $\{i\}$ is the least non-negative
residue class of $i$ modulo $d$. Hence by definition we have
$\{i+kd\}=\{i\}$ for any integer $k\geq0$ and the following equality
holds for any integer $i\in\{1,...,d\}$:
$$\beta_p(i+kd)=\beta_p(i).$$ Then the terms $\frac{1}{i}$ and
$\frac{1}{i+kd}$ of (\ref{Main equivalence}) have the same
coefficient $\beta_p(i)$. For $i\in\{1,...,d\}$ regrouping the
integers $\frac{1}{j}$ such that $j$ lies in the equivalence class
of $i$ modulo $d$ and $j\in\{1,...,\frac{p-1}{2}\}$, the sum in
(\ref{Main equivalence}) can be written
$$\sum_{i=1}^{\frac{p-1}{2}}\beta_p(i)\frac{1}{i}=\sum_{i=1}^{d}\beta_p(i)\alpha_p(i).$$
Then the field $\mathbf{Q}(\sqrt{d})$ is $p$-rational if and only if
$\sum_{i=1}^{d}\beta_p(i)\alpha_p(i)\not\equiv 0\pmod{p}$.\hfill
$\blacksquare$\vskip 6pt

As a consequence we have the following characterization of the
$p$-rationality of the field $\mathbf{Q}(\sqrt{5})$.

\begin{coro}
For every prime $p\equiv1\pmod{5}$, the field $\mathbf{Q}(\sqrt{5})$
is $p$-rational if and only if
\begin{equation}
\alpha_p(1)+\alpha_p(2)-\alpha_p(4)+2\alpha_p(5)\not\equiv0\pmod{p}.
\end{equation}
\end{coro}

\noindent{\textbf{Proof}.} Let $\ell$ be a prime number, then
$(\frac{5}{\ell})=1$ if and only if $\ell\equiv1,4\pmod{5}$, and
$(\frac{5}{\ell})=-1$ if and only if $\ell\equiv2,3\pmod{5}$. Since
$p\equiv1\pmod{5}$, we have for $i\in\{1,...,5\}$,
\begin{equation*}
\beta_p(i)=\sum_{j=1}^{i-1}(\frac{5}{j}),
\end{equation*}
such that $\beta_p(1)=1$, $\beta_p(2)=1$, $\beta_p(3)=0$,
$\beta_p(4)=-1$ and $\beta_p(5)=2$.\hfill $\blacksquare$\vskip 6pt

If we fix the prime number $p$, we obtain a description of the set
of fundamental discriminants $d$ for which the field
$\mathbf{Q}(\sqrt{d})$ is $p$-rational. For the particular cases
$p=3$ and $p=5$ we obtain the following proposition.

\begin{pro}\label{3,5-rational}
Let $d$ be a fundamental discriminant such that $3,5\nmid
(\varepsilon_d-\overline{\varepsilon}_d)^{2}$ then we have the
equivalence:

\begin{enumerate}
    \item $\mathbf{Q}(\sqrt{d})$\;\hbox{is
          3-rational}\;$\Leftrightarrow$\;$\beta_3(1)\not\equiv0\pmod{3}$,
    \item $\mathbf{Q}(\sqrt{d})$\;\hbox{is
          5-rational}\;$\Leftrightarrow$\;$\beta_5(1)\not\equiv2\beta_5(2)\pmod{5}$.
\end{enumerate}
\end{pro}

\noindent{\textbf{Proof}.} By Theorem \ref{Williams}, we have for
$p=3$ the equality
$$\sum_{i=1}^{\frac{3-1}{2}}\beta_3(i)\frac{1}{i}=\beta_3(1),$$
and for $=5$,
$$\sum_{i=1}^{\frac{5-1}{2}}\beta_5(i)\frac{1}{i}=\beta_5(1)+3\beta_5(2).$$
Then the equivalences in $(1)$ and $(2)$ follow from Theorem
\ref{sum of reciprocals}.\hfill $\blacksquare$\vskip 6pt

In general, given an odd prime number $p$, it is not known wether
there exist infinitely many real quadratic fields which are
$p$-rational. This is known for the cases of $p=3$ which is proved
by Dongho Byeon in \cite[Theorem 1]{Byeon}, and the other case is
$p=5$ (see \cite{Assim-Bouazzaoui}). Both cases are proved using
divisibility
properties of Fourier coefficients of half-integer weight modular forms.\\

\noindent{\textbf{Acknowledgement}.} I would like to thank my
advisor J.Assim for his guidance and patience during the preparation
of this paper. Many thanks goes to H.Cohen and B.Abombert for their
help on pariGP computations during the Atelier pariGP in
Besan\c{c}on.

\end{document}